# Solving Boundary Value Problem for a Nonlinear Stationary Controllable System with Synthesizing Control

**Alexander N. Kvitko, Oksana S. Firyulina and Alexey S. Eremin**

*Faculty of Applied Mathematics and Control Processes, Department of Information Systems, Saint-Petersburg State University, Universitetskii prospekt 35, Petergof, Saint Petersburg, Russia, 198504*

An algorithm for constructing a control function that transfers a wide class of stationary nonlinear systems of ordinary differential equations from an initial state to a final state under certain control restrictions is proposed. The algorithm is designed to be convenient for numerical implementation. A constructive criterion of the desired transfer possibility is presented. The problem of an interorbital flight is considered as a test example and it is simulated numerically with the presented method.

## 1. Introduction

One of the problems of mathematical control theory is developing of exact or approximate methods to construct control functions and corresponding trajectories, which connect given points in the phase space. A large amount of publications is devoted to researches in this field, for instance [1–12]. Today boundary value problems (BVPs) are quite well studied for linear and non-linear controllable systems of the special form. However the theory of BVPs for general non-linear controllable systems has not yet been sufficiently developed. The main goal of the authors was to construct an algorithm of solving BVPs for a larger class of non-linear controllable systems of ordinary differential equations in the class of synthesizing controls, which would be numerically stable and easy to implement with computer, and to find a constructive sufficient condition of the solution existence for such problems. This goal was reached by reducing the original problem to a linear non-stationary system of a special form and solving the initial value problem for an auxiliary system of ordinary differential equations. The efficiency of the presented algorithm is demonstrated with numerical simulation of a certain practical problem.

The object of the study is a controllable system of ordinary differential equations (ODEs)

$$\dot{x} = f(x, u), \quad (1.1)$$

where $x = (x_1, \ldots, x_n)^T$ is a vector of length $n$ and $u$ is a vector of same or lesser dimension: $u = (u_1, \ldots, u_r)^T$ and $r \leq n$. We consider the time to satisfy $t \in [0,1]$. The right-hand side

$$f \in C^{4n}(R^n \times R^r; R^n), f = (f_1, \ldots, f_n)^T, \quad (1.2)$$

$$f(0,0) = 0. \quad (1.3)$$

With the denotations

$$A = \frac{\partial f}{\partial x}(0,0), B = \frac{\partial f}{\partial u}(0,0)$$

$$\operatorname{rank} S = n, S = (B, AB, A^2B, \ldots, A^{n-1}B). \quad (1.4)$$

We also consider

$$\|u\| < N. \quad (1.5)$$

**Problem:** To find a pair of functions $x(t) \in C[0,1]$ and $u(t) \in C[0,1]$ that satisfy (1.1) and the conditions

$$x(0) = 0 \text{ and } x(1) = \bar{x}, \bar{x} = (\bar{x}_1, \ldots, \bar{x}_n)^T. \quad (1.6)$$

We say that such pair $x(t), u(t)$ is a solution of the problem (1.1), (1.6).

**Theorem.** Let the conditions (1.2), (1.3) and (1.4) to be satisfied for the right-hand side of (1.1). Then $\exists \varepsilon > 0$ such that $\forall \bar{x} \in R^n$: $\|\bar{x}\| < \varepsilon$ there exists a solution of the problem (1.1), (1.6), which can be found after solving, first, a problem of stabilizing a linear non-stationary system with exponential coefficients and, second, an initial value problem for an auxiliary ODE system.

The main idea of the proof is to use successive changes of independent and dependent variables to reduce the process of solving the original system to the problem of stabilizing a non-linear auxiliary system of ODEs of the special form under constant perturbations. To solve the latter we find a synthesizing control, which provides exponential decrease of the linear auxiliary system fundamental matrix. At the final stage we return to the original variables.

## 2. Auxiliary system construction

We find the function $x(t)$ being a part of the solution of (1.1), (1.6) in the from

$$x_i(t) = a_i(t) + \bar{x}_i, i = 1, \ldots, n. \quad (2.1)$$

In the new variables the system (1.1) and the boundary conditions are written as

$$\dot{a} = f(\bar{x} + a, u), \quad (2.2)$$

$$a(0) = -\bar{x}, a(1) = 0. \quad (2.3)$$

We call a pair of functions $a(t)$ and $u(t)$, which satisfy the system (2.2) and the conditions (2.3), a solution of (2.2), (2.3). Consider a problem to find $a(t) \in C[0,1]$ and $u(t) \in C[0,1]$, which satisfy (2.2), such that

$$a(0) = -\bar{x} \text{ and } a(t) \to 0 \text{ as } t \to 1. \quad (2.4)$$

**Remark 1.** The limit with $t \to 1$ of the solution of (2.2), (2.4) is the solution of (2.2), (2.3).

We make a change of the independent variable $t$ in the system (2.2):

$$t = 1 - e^{-\alpha\tau}, \tau \in [0, +\infty], \quad (2.5)$$

where $\alpha > 0$ is a certain constant value to be determined. Then in terms of $\tau$ (2.2) and (2.4) take form:

$$\frac{dc}{d\tau} = \alpha e^{-\alpha\tau} f(\bar{x} + c, d), \tau \in [0, +\infty], \quad (2.6)$$

$$c(0) = -\bar{x}, c(\tau) \to 0 \text{ as } \tau \to \infty, \quad (2.7)$$

$$c(\tau) = a(t(\tau)), d(\tau) = u(t(\tau)), \\ c = (c_1, \ldots, c_n)^T, d = (d_1, \ldots, d_r)^T. \quad (2.8)$$

We call a pair of functions $c(\tau)$ and $d(\tau)$, which satisfy the system (2.6) with the conditions (2.7), a solution of the problem (2.6), (2.7). With the solution of (2.6), (2.7) one can restore the solution of (2.2), (2.4) with (2.5) and (2.8).

Let's denote

$$\tilde{c} = \bar{x} + \theta_i c, \tilde{d} = \theta_i d, \theta \in [0, 1], i = 1, \ldots, n,$$
$$|k| = \sum_{i=1}^{n} k_i, \quad |m| = \sum_{i=1}^{r} m_i,$$
$$k! = k_1! \cdot \ldots \cdot k_n!, \quad m! = m_1! \cdot \ldots \cdot m_r!$$

Using the property (1.2) and Taylor series expansion of the right-hand side of (1.1) about $(\bar{x}, 0)$, we can rewrite the system (2.6) as:

$$\frac{dc_i}{d\tau} = \alpha e^{-\alpha\tau} f_i(\bar{x}, 0) + \alpha e^{-\alpha\tau} \sum_{j=1}^{n} \frac{\partial f_i}{\partial x_j}(\bar{x}, 0) c_j + \alpha e^{-\alpha\tau} \sum_{j=1}^{r} \frac{\partial f_i}{\partial u_j}(\bar{x}, 0) d_j$$

$$+ \frac{1}{2} \alpha e^{-\alpha\tau} \left[ \sum_{j=1}^{n} \sum_{k=1}^{n} \frac{\partial^2 f_i}{\partial x_j \partial x_k}(\bar{x}, 0) c_j c_k + 2 \sum_{j=1}^{n} \sum_{k=1}^{r} \frac{\partial^2 f_i}{\partial x_j \partial u_k}(\bar{x}, 0) c_j d_k + \sum_{j=1}^{r} \sum_{k=1}^{r} \frac{\partial^2 f_i}{\partial u_j \partial u_k}(\bar{x}, 0) d_j d_k \right]$$

$$+ \alpha e^{-\alpha\tau} \sum_{|k|+|m|=4n-2} \frac{1}{k! m!} \frac{\partial^{|k|+|m|} f_i}{\partial x_1^{k_1} \ldots \partial x_n^{k_n} \partial u_1^{m_1} \ldots \partial u_r^{m_r}}(\bar{x}, 0) c_1^{k_1} \times \ldots \times c_n^{k_n} d_1^{m_1} \times \ldots \times d_r^{m_r} \quad (2.9)$$

$$+ \alpha e^{-\alpha\tau} \sum_{|k|+|m|=4n-1} \frac{1}{k! m!} \frac{\partial^{|k|+|m|} f_i}{\partial x_1^{k_1} \ldots \partial x_n^{k_n} \partial u_1^{m_1} \ldots \partial u_r^{m_r}}(\tilde{c}, \tilde{d}) c_1^{k_1} \times \ldots \times c_n^{k_n} d_1^{m_1} \times \ldots \times d_r^{m_r}, i = 1, \ldots, n.$$

Let's bound the range of $c(\tau)$ with

$$\|c(\tau)\| < C_1, \tau \in [0, \infty). \quad (2.10)$$

We will now shift the functions $c_i(\tau), i = 1, \ldots, n$, several times. Our aim is to get an equivalent system where all the terms in the right-hand side, which do not contain powers of $c$ or $d$ in explicit form, would be of the order $O(e^{-4n\alpha\tau}\|\bar{x}\|)$ as $\tau \to \infty$ and $\|\bar{x}\| \to 0$ in the area (1.5), (2.10).

At the first stage, we change $c_i(\tau)$ to $c_i^{(1)}(\tau)$ by the following rule

$$c_i(\tau) = c_i^{(1)} - e^{-\alpha\tau} f_i(\bar{x}, 0), i = 1, \ldots, n. \quad (2.11)$$

Let $D^{|k|+|m|} f_i \equiv \frac{\partial^{|k|+|m|} f_i}{\partial x_1^{k_1} \ldots \partial x_n^{k_n} \partial u_1^{m_1} \ldots \partial u_r^{m_r}}$ for all $i = 1, \ldots, n$. After substation of (2.11) into the left- and right-hand sized of (2.9) we obtain the equivalent system

$$\frac{dc_i^{(1)}}{d\tau} = -\alpha e^{-2\alpha\tau} \sum_{j=1}^{n} \frac{\partial f_i}{\partial x_j}(\bar{x}, 0) f_j(\bar{x}, 0) + \frac{1}{2} \alpha e^{-3\alpha\tau} \sum_{j=1}^{n} \sum_{k=1}^{n} \frac{\partial^2 f_i}{\partial x_j \partial x_k}(\bar{x}, 0) f_j(\bar{x}, 0) f_k(\bar{x}, 0)$$

$$+ \alpha \left[ e^{-\alpha\tau} \sum_{j=1}^{n} \frac{\partial f_i}{\partial x_j}(\bar{x}, 0) c_j^{(1)} + e^{-2\alpha\tau} \sum_{j=1}^{n} \sum_{k=1}^{n} \frac{\partial^2 f_i}{\partial x_j \partial x_k}(\bar{x}, 0) f_k(\bar{x}, 0) c_j^{(1)} \right]$$

$$+ \alpha \left[ e^{-\alpha\tau} \sum_{j=1}^{r} \frac{\partial f_i}{\partial u_j}(\bar{x}, 0) d_k + e^{-2\alpha\tau} \sum_{j=1}^{n} \sum_{k=1}^{r} \frac{\partial^2 f_i}{\partial x_j \partial u_k}(\bar{x}, 0) f_j(\bar{x}, 0) d_k \right]$$

$$+ \frac{1}{2} \alpha e^{-\alpha\tau} \sum_{j=1}^{n} \sum_{k=1}^{n} \frac{\partial^2 f_i}{\partial x_j \partial x_k}(\bar{x}, 0) c_j^{(1)} c_k^{(1)} + \alpha e^{-\alpha\tau} \sum_{j=1}^{n} \sum_{k=1}^{r} \frac{\partial^2 f_i}{\partial x_j \partial u_k}(\bar{x}, 0) d_k c_j^{(1)} \quad (2.12)$$

$$+ \frac{1}{2} \alpha e^{-\alpha\tau} \sum_{j=1}^{r} \sum_{k=1}^{r} \frac{\partial^2 f_i}{\partial u_j \partial u_k}(\bar{x}, 0) d_j d_k + \cdots$$

$$+ \alpha e^{-\alpha\tau} \sum_{|k|+|m|=4n-2} \frac{1}{k! m!} D^{|k|+|m|} f_i(\bar{x}, 0) \left( c_1^{(1)} - e^{-\alpha\tau} f_1(\bar{x}, 0) \right)^{k_1} \times \ldots$$

$$\times \left( c_n^{(1)} - e^{-\alpha\tau} f_n(\bar{x}, 0) \right)^{k_n} d_1^{m_1} \times \ldots \times d_r^{m_r}$$

$$+ \alpha e^{-\alpha\tau} \sum_{|k|+|m|=4n-1} \frac{1}{k! m!} D^{|k|+|m|} f_i(\tilde{c}, \tilde{d}) \left( c_1^{(1)} - e^{-\alpha\tau} f_1(\bar{x}, 0) \right)^{k_1} \times \ldots$$

$$\times \left(c_n^{(1)} - e^{-\alpha\tau}f_n(\bar{x},0)\right)^{k_n} d_1^{m_1} \times \ldots \times d_r^{m_r}, \; i = 1,\ldots,n.$$

It follows from (2.7) and (2.11) that

$$c_i^{(1)}(0) = -\bar{x}_i + f_i(\bar{x},0), \; i = 1,\ldots,n. \quad (2.13)$$

It is easy to see that in the right-hand side of (2.12) the terms, which do not contain the powers of the components of $c$ or $d$ in explicit form, are bounded with $O(e^{-2\alpha\tau}\|\bar{x}\|)$ as $\tau \to \infty$ and $\|\bar{x}\| \to 0$ in (1.5), (2.10).

At the second stage, we make a change of variables

$$c_i^{(1)}(\tau) = c_i^{(2)}(\tau) + e^{-2\alpha\tau}\phi_i^{(2)}(\bar{x}),$$
$$\phi_i^{(2)}(\bar{x}) = \frac{1}{2}\sum_{j=1}^{n}\frac{\partial f_i}{\partial x_j}(\bar{x},0)f_j(\bar{x},0), \quad (2.14)$$
$$\phi_i^{(2)}(0) = 0, \; i = 1,\ldots,n.$$

With respect to these new variables the original system (2.12) and the initial conditions (2.13) are written as:

$$\begin{aligned}\frac{dc_i^{(2)}}{d\tau} &= \alpha\left[\frac{1}{2}e^{-3\alpha\tau}\sum_{j=1}^{n}\sum_{k=1}^{n}\frac{\partial^2 f_i}{\partial x_j \partial x_k}(\bar{x},0)f_j(\bar{x},0)f_k(\bar{x},0) + e^{-3\alpha\tau}\sum_{j=1}^{n}\frac{\partial f_i}{\partial x_j}(\bar{x},0)\phi_j^{(2)}\right.\\
&\left.-e^{-4\alpha\tau}\sum_{j=1}^{n}\sum_{k=1}^{n}\frac{\partial^2 f_i}{\partial x_j \partial x_k}(\bar{x},0)f_k(\bar{x},0)\phi_j^{(2)} + \frac{1}{2}e^{-5\alpha\tau}\sum_{j=1}^{n}\sum_{k=1}^{n}\frac{\partial^2 f_i}{\partial x_j \partial x_k}(\bar{x},0)\phi_j^{(2)}\phi_k^{(2)}\right]\\
&+\alpha\left[e^{-\alpha\tau}\sum_{j=1}^{n}\frac{\partial f_i}{\partial x_j}(\bar{x},0)c_j^{(2)} - e^{-2\alpha\tau}\sum_{j=1}^{n}\sum_{k=1}^{n}\frac{\partial^2 f_i}{\partial x_j \partial x_k}(\bar{x},0)f_k(\bar{x},0)c_j^{(2)}\right.\\
&\left.+e^{-3\alpha\tau}\sum_{j=1}^{n}\sum_{k=1}^{n}\frac{\partial^2 f_i}{\partial x_j \partial x_k}(\bar{x},0)\phi_k^{(2)}c_j^{(2)}\right]\\
&+\alpha\left[e^{-\alpha\tau}\sum_{k=1}^{r}\frac{\partial f_i}{\partial u_{jk}}(\bar{x},0)d_k - e^{-2\alpha\tau}\sum_{j=1}^{n}\sum_{k=1}^{r}\frac{\partial^2 f_i}{\partial x_j \partial u_k}(\bar{x},0)f_j(\bar{x},0)d_k + e^{-3\alpha\tau}\sum_{j=1}^{n}\sum_{k=1}^{r}\frac{\partial^2 f_i}{\partial x_j \partial u_k}(\bar{x},0)\phi_j^{(2)}d_k\right]\\
&+\frac{1}{2}\alpha e^{-\alpha\tau}\left[\sum_{j=1}^{n}\sum_{k=1}^{n}\frac{\partial^2 f_i}{\partial x_j \partial x_k}(\bar{x},0)c_j^{(2)}c_k^{(2)} + \sum_{j=1}^{n}\sum_{k=1}^{r}\frac{\partial^2 f_i}{\partial x_j \partial u_k}(\bar{x},0)d_k c_j^{(2)} + \sum_{j=1}^{r}\sum_{k=1}^{r}\frac{\partial^2 f_i}{\partial u_j \partial u}(\bar{x},0)d_j d_k\right]\\
&+\alpha e^{-\alpha\tau}\sum_{|k|+|m|=4n-2}\frac{1}{k!m!}D^{|k|+|m|}f_i(\bar{x},0)\left(c_1^{(2)} + e^{-2\alpha\tau}\phi_1^{(2)} - e^{-\alpha\tau}f_1(\bar{x},0)\right)^{k_1}\times\ldots\\
&\times\left(c_n^{(2)} + e^{-2\alpha\tau}\phi_n^{(2)} - e^{-\alpha\tau}f_n(\bar{x},0)\right)^{k_n} d_1^{m_1} \times\ldots\times d_r^{m_r}\\
&+\alpha e^{-\alpha\tau}\sum_{|k|+|m|=4n-1}\frac{1}{k!m!}D^{|k|+|m|}f_i(\tilde{c},\tilde{d})\left(c_1^{(2)} + e^{-2\alpha\tau}\phi_1^{(2)} - e^{-\alpha\tau}f_1(\bar{x},0)\right)^{k_1}\times\ldots\\
&\times\left(c_n^{(2)} + e^{-2\alpha\tau}\phi_n^{(2)} - e^{-\alpha\tau}f_n(\bar{x},0)\right)^{k_n} d_1^{m_1}\times\ldots\times d_r^{m_r}, \; i=1,\ldots,n.\end{aligned} \quad (2.15)$$

$$c_1^{(2)}(0) = -\bar{x}_i + f_i(\bar{x},0) - \phi_i^{(2)}(\bar{x}), \; i = 1,\ldots,n. \quad (2.16)$$

Comparing to the previous variables shift, the right-hand side terms of (2.15), which do not contain the powers of the components of $c$ or $d$ in explicit form, are now of the order $O(e^{-3\alpha\tau}\|\bar{x}\|)$ as $\tau \to \infty$ and $\|\bar{x}\| \to 0$ in the area (1.5), (2.10).

By induction, at the $k$-th stage, using (2.11)–(2.16) we have the necessary shift of the form:

$$c_i^{(k-1)}(\tau) = c_i^{(k)} + e^{-k\alpha\tau}\phi_i^{(k)}(\bar{x}),$$
$$\phi_i^{(k)}(0) = 0, \; i = 1,\ldots,n. \quad (2.17)$$

We apply (2.17) $4n-1$ times and collect the terms, which are linear in respect to the components of $c^{(4n-1)}$ and include the coefficients $e^{-i\alpha\tau}, \; i=1,\ldots,n$, and also the terms, which are linear in respect to the components of $d$ and include the coefficients $e^{-i\alpha\tau}, \; i=1,\ldots,2n$. Now we have the system, which according to (2.12)–(2.17) can be written in vector form as:

$$\begin{aligned}\frac{dc^{(4n-1)}}{d\tau} &= Pc^{(4n-1)} + Qd\\
&+R_1(c^{(4n-1)},d,\bar{x},\tau) + R_2(c^{(4n-1)},d,\bar{x},\tau)\\
&+R_3(c^{(4n-1)},d,\tau) + R_4(\bar{x},c^{(4n-1)},d,\tau),\\
R_1 &= (R_1^1,\ldots,R_1^n)^T, \; R_2 = (R_2^1,\ldots,R_2^n)^T,\\
R_3 &= (R_3^1,\ldots,R_3^n)^T, \; R_4 = (R_4^1,\ldots,R_4^n)^T.\end{aligned} \quad (2.18)$$

The functions $R_1^i$ consist all the terms which are linear in respect to the components of $c^{(4n-1)}$ with coefficients $e^{-i\alpha\tau}, \; i \geq n+1$, and also the terms of the last sum of the right-hand side, for which $|m|=0$ and $|k|=1$. The functions $R_2^i$ consist all the terms which are linear in respect to the components of $d$ with coefficients $e^{-i\alpha\tau}, \; i \geq 2n+1$, and also the terms of the last sum of the right-hand side, for which $|m|=1$ and $|k|=0$. In $R_3^i$ all the terms, which are non-linear in respect of the components of $c^{(4n-1)}$ or $d$, are contained. Finally, the functions $R_4^i$ include all the terms,

which do not have powers of $c^{(4n-1)}$ and $d$ components. The functions $P$ and $Q$ have form

$$P(\bar{x}) = \alpha e^{-\alpha\tau}\left(P_1(\bar{x}) + e^{-\alpha\tau}P_2(\bar{x}) + \cdots + e^{-(n-1)\alpha\tau}P_{n-1}(\bar{x})\right),$$
$$P_1(\bar{x}) = \frac{\partial f}{\partial x}(\bar{x}, 0),\ P_1(0) = A,$$
$$Q(\bar{x}) = \alpha e^{-\alpha\tau}\left(Q_1(\bar{x}) + e^{-\alpha\tau}Q_2(\bar{x}) + \cdots + e^{-(2n-1)\alpha\tau}Q_{2n-1}(\bar{x})\right),$$
$$Q_1(\bar{x}) = \frac{\partial f}{\partial u}(\bar{x}, 0),\ Q_1(0) = B.$$
(2.19)

$$c^{(4n-1)}(0) = -\bar{x} + f(\bar{x}) - \phi^{(2)}(\bar{x}) - \cdots - \phi^{(n-1)}(\bar{x}),$$
$$\phi^{(i)} = (\phi_1^{(i)}, \ldots, \phi_n^{(i)})^T,\ i = 1, \ldots, 4n-1,$$
$$\phi^{(i)}(0) = 0.$$
(2.20)

## 3. Right hand side terms evaluation

It follows from the construction of (2.18) that in the area (1.5), (2.10) the following estimations are true

$$\|P_i(\bar{x})\| \to 0, \|Q_j(\bar{x})\| \to 0 \text{ as } \|\bar{x}\| \to 0, i = 2, \ldots, n-1,\ j = 2, \ldots, 2n-1;$$
(3.1)

$$\|R_1(c^{(4n-1)}, d, \bar{x}, \tau)\| \le e^{-(n+1)\alpha\tau}L_1(\bar{x})\|c^{(4n-1)}\|,$$
$$\|R_2(c^{(4n-1)}, d, \bar{x}, \tau)\| \le e^{-(2n+1)\alpha\tau}L_2(\bar{x})\|d\|,\ L_1 > 0, L_2 > 0;$$
(3.2)

$$\|R_3(c^{(4n-1)}, d, \tau)\| \le e^{-\alpha\tau}L_3\left(\|c^{(4n-1)}\|^2 + \|d\|^2\right),$$
$$L_3 > 0.$$
(3.3)

Moreover, the conditions (1.2) and (1.3) lead to

$$\|R_4(c^{(4n-1)}, d, \bar{x}, \tau)\| \le L_4\|\bar{x}\|e^{-4n\alpha\tau},\ L_4 > 0. \quad (3.4)$$

The estimation (3.4) follows from the representation

$$f(\bar{x}, 0) = \frac{\partial f}{\partial x}(\theta\bar{x}, 0)\bar{x},\ \theta = (\theta_1, \ldots, \theta_n)^T,$$
$$\theta_i \in [0,1],\ i = 1, \ldots, n.$$

**Remark 2.** Let's denote the $i$-th column of the matrix $Q_1$ as $q_1^i$. Construct a matrix

$$S_1 = \{q_1^1, P_1 q_1^1, \ldots, P_1^{k_1-1} q_1^1, q_1^2, P_1 q_1^2, \ldots, P_1^{k_2-1} q_1^2, \ldots, q_1^r, \ldots, P_1^{k_r-1} q_1^r\}.$$

Here for each $j = 1, \ldots, r$ $k_j$ is the maximal number of columns of the form $q_1^j, \ldots, P_1^{k_j-1} q_1^j$ such that all the vectors $q_1^1, P_1 q_1^1, \ldots, P_1^{k_1-1} q_1^1, q_1^2, P_1 q_1^2, \ldots, P_1^{k_2-1} q_1^2, q_1^r, \ldots, P_1^{k_r-1} q_1^r$ are linearly independent.

From the conditions (1.4) and (2.19) it follows, that there exists $\varepsilon_1 > 0$ so that $\text{rank}\, S_1 = n$ $\forall \bar{x} \in R^n$ that satisfies $\|\bar{x}\| < \varepsilon_1$.

We now consider the system

$$\frac{dc^{(4n-1)}}{d\tau} = Pc^{(4n-1)} + Qd. \quad (3.5)$$

## 4. Auxiliary lemma

**Lemma.** Let the conditions (1.2) and (1.4) be satisfied for the system (1.1). Then $\exists \bar{\varepsilon} > 0,\ \bar{\varepsilon} < \varepsilon_1$ such that $\forall \bar{x} \in R^n$: $\|\bar{x}\| < \bar{\varepsilon}$ there exists a control $d(\tau)$ of the form

$$d(\tau) = M(\tau)c^{(4n-1)} \quad (4.1)$$

that provides an exponential decrease of the fundamental matrix in (3.5).

*Proof.* We use Krasosvsky linear non-stationary systems stabilization method in the proof. Let $L_1^j, j = 1, \ldots, r$, to be the $j$-th column of $Q$. Construct a matrix

$$S_2 = \{L_1^1, L_2^1, \ldots, L_{k_1}^1, L_1^2, \ldots, L_{k_2}^2, \ldots, L_1^r, \ldots, L_{k_r}^r\},$$
$$L_i^j = PL_{i-1}^j - \frac{dL_{i-1}^j}{d\tau},\ j = 1, \ldots, r,\ i = 2, \ldots, k_j.$$
(4.2)

Here for each $j$ $k_j$ is the maximal number of columns $L_1^j$, …, $L_{k_j}^j$ such that the vectors $L_1^1, L_2^1, \ldots, L_{k_1}^1, L_1^2, \ldots, L_{k_2}^2, L_1^r$, …, $L_{k_r}^r$ are linearly independent. Let's show that

$$\text{rank}\, S_2 = n. \quad (4.3)$$

Let $\bar{L}_1^j, j = 1, \ldots, r$, to be the $j$-th column of the matrix $\alpha e^{-\alpha\tau}Q_1$. Consider the matrix

$$S_3 = \{\bar{L}_1^1, \bar{L}_2^1, \ldots, \bar{L}_{k_1}^1, \bar{L}_1^2, \ldots, \bar{L}_{k_2}^2, \ldots, \bar{L}_1^r, \ldots, \bar{L}_{k_r}^r\},$$
$$\bar{L}_i^j = \alpha e^{-\alpha\tau} P_1 \bar{L}_{i-1}^j - \frac{d\bar{L}_{i-1}^j}{d\tau},\ j = 1, \ldots, r,\ i = 2, \ldots, k_j$$
(4.4)

where $k_j$ are determined the same way as for $S_2$. The conditions (1.2), (2.19) and (3.1) provide that $S_2 \to S_3$ as $\|\bar{x}\| \to 0$. This lead to the existence of $\varepsilon_2 > 0$: $\varepsilon_2 < \varepsilon_1$ such that $\forall \bar{x} \in R^n$: $\|\bar{x}\| < \varepsilon_2$

$$\text{rank}\, S_2 = \text{rank}\, S_3. \quad (4.5)$$

Now, taking in account the Remark 2 and the equality (4.5), we can prove by contradiction that $\forall \bar{x} \in R^n$: $\|\bar{x}\| < \varepsilon_2$

$$\text{rank}\, S_3 = n. \quad (4.6)$$

From (4.5) and (4.6) it follows that the condition (4.3) is satisfied in the area $\|\bar{x}\| < \varepsilon_2$. Moreover, in that area the structure of $S_3$ provides the estimation

$$\|S_2^{-1}\| = O(e^{n\alpha\tau}),\ \tau \to \infty. \quad (4.7)$$

With use of (4.3), change the variable $c^{(4n-1)}$ according to the expression

$$c^{(4n-1)} = S_2(\tau)y. \quad (4.8)$$

As a result we get the system

$$\frac{dy}{d\tau} = S_2^{-1}\left(PS_2 - \frac{dS_2}{d\tau}\right)y + S_2^{-1}Qd. \quad (4.9)$$

According to [13] for the first term of the right-hand side of (4.9) we have

$$S_2^{-1}\left(PS_2 - \frac{dS_2}{d\tau}\right) = \{\bar{e}_2, \ldots, \bar{e}_{k_1}, \bar{\varphi}_{k_1}(\tau), \ldots,$$
$$\bar{e}_{k_1+\cdots+k_{r-1}+2}, \ldots, \bar{e}_{k_1+\cdots+k_r}, \bar{\varphi}_{k_r}(\tau)\}.$$

where $\bar{e}_i$ is a zero vector of size $n$ with the only unit at $i$-th place,

$$\bar{\varphi}_{k_j} = \left(-\varphi_{k_1}^1, \ldots, -\varphi_{k_1}^{k_1}, \ldots, -\varphi_{k_j}^1, \ldots, -\varphi_{k_j}^{k_j}, 0, \ldots, 0\right)_{n\times 1}^T$$

and $\varphi_{k_j}^i$ are the coefficient of $L_{k_j+1}^j$ expansion into the sum of the vectors $L_1^1, L_2^1, \ldots, L_{k_1}^1, L_1^2, \ldots, L_{k_2}^2, L_1^r, \ldots, L_{k_r}^r$, (notice that $\sum_{j=1}^r k_j = n$) i.e.

$$L_{k_j+1}^j = -\sum_{i=1}^{k_1}\varphi_{k_1}^i(\tau)L_i^1 - \cdots - \sum_{i=1}^{k_j}\varphi_{k_j}^i(\tau)L_i^j. \quad (4.10)$$

The second term is

$$S_2^{-1}Q = \{\bar{e}_1, \ldots, \bar{e}_{k_i+1}, \ldots, \bar{e}_{\gamma+1}\}, \gamma = \sum_{i=1}^{r-1}k_i.$$

Consider the stabilization of the system

$$\frac{dy_{k_i}}{d\tau} = \left\{\bar{e}_2^{k_i}, \ldots, \bar{e}_{k_i}^{k_i}, \bar{\bar{\varphi}}_{k_i}\right\}y_{k_i} + \bar{e}_1^{k_i}d_i,$$
$$y_{k_i} = \left(y_{k_i}^1, \ldots, y_{k_i}^{k_i}\right)_{k_i \times 1}^T, \quad i = 1, \ldots, r,$$
(4.11)

where $\bar{e}_i^{k_i}$ is a zero vector of size $k_i$ with the only unit at $i$-th place, and $\bar{\bar{\varphi}}_{k_i} = \left(-\varphi_{k_i}^1, \ldots, -\varphi_{k_i}^{k_i}\right)_{k_i\times 1}^T$.

Let $y_{k_i}^{k_i} = \psi$. The phase variables of (4.11) are connected to $\psi(\tau)$ and its derivatives with the equalities:

$$y_{k_i}^{k_i} = \psi,$$
$$y_{k_i}^{k_i-1} = \psi^{(1)} + \varphi_{k_i}^{k_i}\psi,$$
$$y_{k_i}^{k_i-2} = \psi^{(2)} + \varphi_{k_i}^{k_i}\psi^{(1)} + \left(\frac{d\varphi_{k_i}^{k_i}}{d\tau} + \varphi_{k_i}^{k_i-1}\right)\psi, \quad (4.12)$$
$$\ldots$$
$$y_{k_i}^1 = \psi^{(k_i-1)} + r_{k_i-2}(\tau)\psi^{(k_i-2)} + \cdots + r_1(\tau)\psi^{(1)} + r_0(\tau)\psi$$

The differentiation of the last equality in (4.12) together with (4.11) leads to

$$\psi^{(k_i)} + \varepsilon_{k_i-1}(\tau)\psi^{(k_i-1)} + \ldots + \varepsilon_0(\tau)\psi = d_i, \quad (4.13)$$
$$i = 1, \ldots, r.$$

The functions $r_{k_i-2}(\tau), \ldots, r_0(\tau), \varepsilon_{k_i-1}(\tau), \ldots, \varepsilon_0(\tau)$ in (4.12) and (4.13) are linear combinations of the functions $\varphi_{k_j}^i(\tau), i = 1, \ldots, k_j, j = 1, \ldots, r$ and their derivatives.

**Remark 3.** From the structure of the matrices $P$ and $Q$ it follows (see (2.19) and the representation (4.10)) that the functions $\varphi_{k_i}^{k_i}(\tau), \ldots, \varphi_{k_i}^1(\tau)$ and their derivatives, as well as $r_{k_i-2}(\tau), \ldots, r_0(\tau), \varepsilon_{k_i-1}(\tau), \ldots, \varepsilon_0(\tau)$ are bounded. Other elements of the columns $\bar{\varphi}_{k_j}, j = 1, \ldots, r,$ obey the estimation $O(e^{(n-1)\alpha\tau}), \tau \to \infty$.

Let

$$d_i = \sum_{j=1}^{k_i}\left(\varepsilon_{k_i-j}(\tau) - \gamma_{k_i-j}\right)\psi^{(k_i-j)}, \quad (4.14)$$
$$i = 1, \ldots, r,$$

where $\gamma_{k_i-j}, j = 1, \ldots, k_i$ are chosen so that the $\lambda_{k_i}^1, \ldots, \lambda_{k_i}^{k_i}$ of the equations

$$\lambda^{k_i} + \gamma_{k_i-1}\lambda^{k_i-1} + \cdots + \gamma_0 = 0, i = 1, \ldots, r,$$

to satisfy the conditions

$$\lambda_{k_i}^i \neq \lambda_{k_i}^j \text{ if } i \neq j \text{ and } \lambda_{k_i}^i < -(2n+1)\alpha - 1 \quad (4.15)$$
$$\text{for } j = 1, \ldots k_i, i = 1, \ldots r.$$

Due to (4.12), (4.15) and the Remark 3, the control $d_i = \delta_{k_i}T_{k_i}^{-1}y_{k_i}, \ i = 1, \ldots, r$ provides the exponential decrease of the system (4.9) solutions. Switching back to the original variables in (4.14) and using (4.8), we have

$$d_i = \delta_{k_i}T_{k_i}^{-1}S_{2k_i}^{-1}c^{(4n-1)}, i = 1, \ldots, r, \quad (4.16)$$

where $\delta_{k_i} = (\varepsilon_{k_i-1}(\tau) - \gamma_{k_i-1}, \ldots, \varepsilon_0(\tau) - \gamma_0); T_{k_i}$ is the inequality (4.12) matrix, which means that $y_{k_i} = T_{k_i}\bar{\psi}, \bar{\psi} = \left(\psi^{(k_i-1)}, \ldots, \psi\right)^T; S_{2k_i}^{-1}$ is the matrix which consists of the corresponding $k_i$ rows of $S_2^{-1}$. The found control can be written in the form (4.1), where

$$M(\tau) = \delta_k T_k^{-1}S_{2k}^{-1} \equiv \left(\delta_{k_1}T_{k_1}^{-1}S_{2k_1}^{-1}, \ldots, \delta_{k_r}T_{k_r}^{-1}S_{2k_r}^{-1}\right)^T.$$

Denote the fundamental matrix of the system (4.13) closed with the control (4.14) via $\Psi(\tau), (\Psi(0) = I,$ an identity matrix). It is obvious that the elements of $\Psi(\tau)$ are exponential functions of negative argument or their derivatives.

Consider the system (3.15) with the control (4.16)

$$\frac{dc^{(4n-1)}}{d\tau} = D(\tau)c^{(4n-1)}, D(\tau) = P(\tau) + Q(\tau)M(\tau). \quad (4.17)$$

Introduce a block-diagonal matrix $T(\tau)$. Its diagonal blocks are matrices $T_{k_i}, i = 1, \ldots, r$. Then, corresponding to (4.8) and (4.12), the fundamental matrix $\Phi(\tau), \Phi^{-1}(0) = I$, of the system (4.17) has form

$$\Phi(\tau) = S_2(\tau)T(\tau)\Psi(\tau)T^{-1}(0)S_2^{-1}(0). \quad (4.18)$$

The estimation

$$\|\Phi(\tau)\| \leq Ke^{-\alpha\tau}e^{-\lambda\tau}, \lambda > 0, K > 0 \quad (4.19)$$

follows from (4.18), the structure of matrices $S_2(\tau)$ and $\Psi(\tau)$ and the Remark 3.

Let $\bar{\varepsilon} = \varepsilon_2$. Then the correctness of the Lemma becomes clear from (4.19).

Besides, on base of (4.7), (4.16), (4.18) and the Remark 3 we have

$$\|\Phi(\tau)\Phi^{-1}(t)\| \leq K_1 e^{-\lambda(\tau-t)}e^{(n-1)\alpha t}, \quad (4.20)$$

$t \leq \tau, \tau \in [0, \infty), K > 0$;
$\|M(\tau)\| = O(e^{n\alpha\tau}), \tau \to \infty$.

## 5. Main theorem proof

The system (2.18) with the control (4.16) has form

$$\frac{dc^{(4n-1)}}{d\tau} = D(\tau)c^{(4n-1)}$$
$$+ R_1\big(c^{(4n-1)}, M(\tau)c^{(4n-1)}, \bar{x}, \tau\big)$$
$$+ R_2\big(c^{(4n-1)}, M(\tau)c^{(4n-1)}, \bar{x}, \tau\big) \quad (5.1)$$
$$+ R_3\big(c^{(4n-1)}, M(\tau)c^{(4n-1)}, \tau\big)$$
$$+ R_4\big(\bar{x}, c^{(4n-1)}, M(\tau)c^{(4n-1)}, \tau\big)$$

Change of the variables

$$c^{(4n-1)} = ze^{-3n\alpha\tau}, \quad c^{(4n-1)}(0) = z(0), \quad (5.2)$$

leads to the following representation

$$\frac{dz}{d\tau} = C(\tau)z + e^{3n\alpha\tau}R_1(e^{-3n\alpha\tau}z, M(\tau)e^{-3n\alpha\tau}z, \bar{x}, \tau)$$
$$+ e^{3n\alpha\tau}R_2(e^{-3n\alpha\tau}z, M(\tau)e^{-3n\alpha\tau}z, \bar{x}, \tau) \quad (5.3)$$
$$+ e^{3n\alpha\tau}R_3(e^{-3n\alpha\tau}z, M(\tau)e^{-3n\alpha\tau}z, \tau)$$
$$+ e^{3n\alpha\tau}R_4(\bar{x}, e^{-3n\alpha\tau}z, M(\tau)e^{-3n\alpha\tau}z, \tau).$$
$$C(\tau) = D(\tau) + 3n\alpha E.$$

Let's show that all solutions of (5.3) with initial values (5.2) that start in a sufficiently small neighborhood of zero decrease exponentially. Let $\Phi_1(\tau)$, $\Phi_1^{-1}(0) = I$, be the fundamental matrix of the system $\frac{dz}{d\tau} = C(\tau)z$. Then on base of (4.19), (4.20) and (5.2) we have

$$\|\Phi_1(\tau)\| \leq Ke^{-\beta\tau},$$
$$\|\Phi_1(\tau)\Phi_1^{-1}(t)\| \leq K_1 e^{-\beta(\tau-t)}e^{(n-1)\alpha t}, \quad (5.4)$$
$$\beta = \lambda - 3n\alpha.$$

Let's choose $\alpha$ to provide

$$\beta > 0. \quad (5.5)$$

The solution of (5.3) with initial values (2.20) can be presented as

$$z(\tau) = \Phi_1(\tau)\Phi_1^{-1}(\tau_1)z(\tau_1) +$$
$$\int_{\tau_1}^{\tau} \Phi_1(\tau)\Phi_1^{-1}(t)e^{3n\alpha t} \times$$
$$[R_1(e^{-3n\alpha t}z, M(t)e^{-3n\alpha t}z, \bar{x}, t) + \quad (5.6)$$
$$R_2(e^{-3n\alpha t}z, M(t)e^{-3n\alpha t}z, \bar{x}, t) +$$
$$R_3(e^{-3n\alpha t}z, M(t)e^{-3n\alpha t}z, t) +$$
$$R_4(\bar{x}, e^{-3n\alpha t}z, M(t)e^{-3n\alpha t}z, t)]dt,$$
$$\text{for } \tau \in [\tau_1, \infty),$$

$$z(\tau) = \Phi_1(\tau)c^{(4n-1)}(0) +$$
$$\int_0^{\tau} \Phi_1(\tau)\Phi_1^{-1}(t)e^{3n\alpha t} \times$$
$$[R_1(e^{-3n\alpha t}z, M(t)e^{-3n\alpha t}z, \bar{x}, t) + \quad (5.7)$$
$$R_2(e^{-3n\alpha t}z, M(t)e^{-3n\alpha t}z, \bar{x}, t) +$$
$$R_3(e^{-3n\alpha t}z, M(t)e^{-3n\alpha t}z, t) +$$
$$R_4(\bar{x}, e^{-3n\alpha t}z, M(t)e^{-3n\alpha t}z, t)]dt,$$
$$\text{for } \tau \in [0, \tau_1).$$

Now, the use (3.2), (3.3), (3.4), (4.1), (5.2) and (5.4) after changing $c$ to $z$, gives the following estimations in the area (1.5), (2.10) from (5.6) and (5.7):

$$\|z(\tau)\| \leq Ke^{-\beta\tau}\|\Phi_1^{-1}(\tau_1)z(\tau_1)\|$$
$$+ \int_{\tau_1}^{\tau} e^{-\beta(\tau-t)}K_1(\bar{L}e^{-\alpha t}\|z(t)\| + L_4\|\bar{x}\|e^{-\alpha t})dt, \quad (5.8)$$
$$\text{for } \tau \in [\tau_1, \infty); \text{ and}$$

$$\|z(\tau)\| \leq Ke^{-\beta\tau}\|c^{(4n-1)}(0)\|$$
$$+ \int_0^{\tau} e^{-\beta(\tau-t)}K_1(\bar{L}e^{-\alpha t}\|z(t)\| + L_4\|\bar{x}\|e^{-\alpha t})dt, \quad (5.9)$$
$$\text{for } \tau \in [0, \tau_1), \bar{L} > 0.$$

The constant $\bar{L}$ depends on (1.5), (2.10).

We apply the well-known result from [14] to the inequalities (5.8) and (5.9) and obtain

$$\|z(\tau)\| \leq Ke^{-\mu\tau}\|\Phi^{-1}(\tau_1)z(\tau_1)\|$$
$$+ K_1 \int_{\tau_1}^{\tau} e^{-\mu(\tau-t)}L_4\|\bar{x}\|e^{-\alpha t}dt,$$
$$\text{for } \tau \in [\tau_1, \infty), \mu = \beta - K_1\bar{L}e^{-\alpha\tau_1}; \text{ and} \quad (5.10)$$
$$\|z(\tau)\| \leq Ke^{-\mu_1\tau}\|c^{(4n-1)}(0)\|$$
$$+ K_1 \int_0^{\tau} e^{-\mu_1(\tau-t)}L_4\|\bar{x}\|e^{-\alpha t}dt,$$
$$\text{for } \tau \in [0, \tau_1], \mu_1 = \beta - K_1\bar{L}.$$

Using (5.5), we fix $\tau_1 > 0$ so, that the inequality $\mu > 0$ is satisfied.

Now let's bound the choice of $\alpha > 0$ with the condition $\alpha < \mu$. Then after the integration in the right-hand sides of (5.10) we have

$$\|z(\tau)\| \leq Ke^{-\mu\tau}\|\Phi_1^{-1}(\tau_1)\|\|(z(\tau_1)\| +$$
$$K_2 e^{-\alpha\tau}L_4\|\bar{x}\|, \text{ for } \tau \in [\tau_1, \infty), \text{ and}$$
$$\|z(\tau)\| \leq K_3\|c^{(4n-1)}(0)\| + K_4 L_4\|\bar{x}\|,$$
$$\text{for } \tau \in [0, \tau_1].$$
$$\text{Not that all } K_i > 0.$$

Based on (2.20), (3.4) two last estimations can be written as a single inequality in the area $\|\bar{x}\| < \bar{\varepsilon}$:

$$\|z(\tau)\| \leq K_5 e^{-\alpha\tau}\|\bar{x}\|, \tau \in [0, \infty), K_5 > 0. \quad (5.11)$$

The constant $K_5$ depends on the area (1.5), (2.10).

With use of (4.1), (4.20), (5.2) and (5.11) we estimate $\|d(\tau)\|$:

$$\|d(\tau)\| \leq \|M(\tau)\|e^{-3n\alpha\tau}K_5 e^{-\alpha\tau}\|\bar{x}\| \leq$$
$$K_6 e^{-(2n-1)\alpha\tau}\|\bar{x}\|, K_6 > 0, \text{ for } \tau \in [0, \infty). \quad (5.12)$$

Now we can set the value $\varepsilon$ from the Theorem formulation, to be $\min\left\{\bar{\varepsilon}, \frac{C_1}{K_5}, \frac{N}{K_6}\right\}$. This shows that the solution of the system (5.3) with initial values (5.2), (2.20) doesn't leave the area $\|z\| < C_1$ for $\|\bar{x}\| < \varepsilon$ and decreases exponentially. Besides, by (5.12) the corresponding function $d(\tau)$ obeys the restriction (1.5). Moreover, if we

substitute the solution to (5.2) and (4.1) and return to the variables $c(\tau)$ according to (2.11), (2.14)–(2.17), we get the solution of (2.6), (2.7).

Then we write everything in the original variables with (2.8), (2.5) and finally (2.1). Now, according to the Remark 1, passing to the limit as $t \to 1$ gives the solution of the original problem (1.1), (1.6). This proves the Theorem.

**Remark 4.** Consider the system

$$\dot{x} = f(x, u) + F \qquad (5.13)$$

where $x = (x_1, \ldots, x_n)^T$, $x \in R^n$; $u = (u_1, \ldots, u_r)^T$, $u \in R^r$; $t \in [0, 1]$, $r \leq n$, $F = (F_1, \ldots, F_1)^T$, $F \in R^n$. $F$ is a constant perturbation. From the Theorem follows the

**Corollary.** Let the conditions (1.2)–(1.4) to be satisfied for the system (5.13). Then there exists $\varepsilon > 0$ so that for any $\bar{x} \in R^n$ and $F \in R^n$ if $\|\bar{x}\| < \varepsilon$ and $\|F\| < \varepsilon$ a solution of (5.13), (1.6) exists and it can be found after stabilizing a linear non-stationary system with exponential coefficients and solving the initial value problem for an auxiliary system of ODEs.

## 6. Practical example

To demonstrate the effectiveness of the proposed method we consider a problem of transferring a material point moving in a central gravitational field to a desired circular orbit with jet power.

According to [12] the system in deviations from the prescribed circular motion and the condition (1.6) take form

$$\dot{x}_1 = x_2,$$

$$\dot{x}_2 = -\frac{\nu}{(r_0 + x_1)^2} + \frac{(\chi_0 + x_3)^2}{(r_0 + x_1)^3} + a_r u, \qquad (6.1)$$

$$\dot{x}_3 = (x_1 + r_0) a_\psi u,$$

$$x_i(0) = 0, x_i(1) = \bar{x}_i, i = 1, 2, 3. \qquad (6.2)$$

Here $x_1 = r - r_0$, $x_2 = \dot{r}$, $x_3 = \chi - \chi_0$; $\chi_0 = (\nu r_0)^{\frac{1}{2}}$, $u = \dot{m}/m$; $r_0$ is the circular orbit radius; $\dot{r}$ is the radial velocity; $\chi$ is the generalized momentum; $a_r$ and $a_\psi$ are the relative velocity vector projections onto the radial and tangential directions respectively (they are constant); $m$ is the mass and $\dot{m}$ is its change rate; $\nu = \nu^0 M$ where $\nu^0$ is the universal gravitational constant; $M$ is the mass of the Earth; vector $x = (x_1, x_1, x_3)^T$; and the control $u \in R^1$.

The system (2.6) and conditions (2.7) for the problem have form

$$\frac{dc_1}{d\tau} = e^{-\alpha\tau} c_2,$$

$$\frac{dc_2}{d\tau} = e^{-\alpha\tau} g(c_1, c_3) + e^{-\alpha\tau} a_r d, \qquad (6.3)$$

$$\frac{dc_3}{d\tau} = e^{-\alpha\tau} (c_1 + \bar{x}_1 + r_0) a_\psi d,$$

with

$$g(c_1, c_3) = -\frac{\nu}{(r_0 + c_1 + \bar{x}_1)^2} + \frac{(\chi_0 + c_3 + \bar{x}_3)^2}{(r_0 + c_1 + \bar{x}_1)^3}$$

and

$$c_i(0) = -\bar{x}_i, c_i(\tau) \to 0 \text{ as } \tau \to \infty, i = 1, 2, 3. \qquad (6.4)$$

In the following we assume $\bar{x} = (\bar{x}_1, 0, 0)^T$. It is necessary to make seven transforms of the type (2.17) to solve the problem (6.3), (6.4):

$$c_2 = z_2 - e^{-\alpha\tau} g(\bar{x}_1), c_1 = z_1 + \frac{1}{2} e^{-2\alpha\tau} g(\bar{x}_1),$$

$$z_2 = w_2 - \frac{1}{6} e^{-3\alpha\tau} \frac{\partial g}{\partial c_1}(\bar{x}_1) g(\bar{x}_1),$$

$$z_1 = w_1 + \frac{1}{24} e^{-4\alpha\tau} \frac{\partial g}{\partial c_1}(\bar{x}_1) g(\bar{x}_1),$$

$$w_2 = v_2 - e^{-5\alpha\tau} \bar{g}, w_1 = v_1 + e^{-6\alpha\tau} \bar{g},$$

$$v_2 = u_2 - e^{-7\alpha\tau} \bar{\bar{g}},$$

$$\bar{g} = \frac{1}{120} \left(\frac{\partial g}{\partial c_1}(\bar{x}_1)\right)^2 g(\bar{x}_1) - \frac{1}{40} \frac{\partial^2 g}{\partial c_1^2}(\bar{x}_1) g^2(\bar{x}_1),$$

$$\bar{\bar{g}} = \frac{1}{720} \left(\frac{\partial g}{\partial c_1}(\bar{x}_1)\right)^2 g(\bar{x}_1) + \frac{1}{240} \frac{\partial^2 g}{\partial c_1^2}(\bar{x}_1) g^2(\bar{x}_1), \qquad (6.5)$$

$$\bar{\bar{\bar{g}}} = \frac{1}{7} \left\{\frac{1}{720} \left(\frac{\partial g}{\partial c_1}(\bar{x}_1)\right)^3 g(\bar{x}_1)\right.$$

$$+ \frac{1}{48} \frac{\partial^3 g}{\partial c_1^3}(\bar{x}_1) g^3(\bar{x}_1)$$

$$\left. + \frac{1}{40} \frac{\partial g}{\partial c_1}(\bar{x}_1) \frac{\partial^2 g}{\partial c_1^2}(\bar{x}) g^2(\bar{x}_1)\right\},$$

The matrices $P$ and $Q$ and the system analogous to (3.5), look like

$$\frac{d\bar{c}}{d\tau} = P\bar{c} + Qd, \bar{c} = (v_1, u_2, c_3)^T,$$

$$P = \alpha e^{-\alpha\tau} \begin{Vmatrix} 0 & 1 & 0 \\ a_{21} & 0 & a_{23} \\ 0 & 0 & 0 \end{Vmatrix},$$

$$Q = \begin{Vmatrix} 0 \\ \alpha e^{-\alpha\tau} a_r \\ \alpha e^{-\alpha\tau} b + \frac{\alpha}{2} e^{-3\alpha\tau} a_\psi \end{Vmatrix}, \qquad (6.6)$$

$$b = (r_0 + \bar{x}_1) a_\psi, a_{21} = \frac{\partial g}{\partial c_1}(\bar{x}_1), a_{23} = \frac{\partial g}{\partial c_3}(\bar{x}_1).$$

To stabilize the system (6.6) we construct the matrix $S = \{L_1, L_2, L_3\}$, with $L_1 = Q$, $L_2 = PL_1 - \frac{d}{d\tau} L_1$, $L_3 = PL_2 - \frac{d}{d\tau} L_2$, i.e.

$$L_1 = \begin{pmatrix} 0 \\ \alpha e^{-\alpha\tau} \\ \alpha e^{-\alpha\tau}(b + e^{-2\alpha\tau} g) \end{pmatrix}$$

$$L_2 = \begin{pmatrix} \alpha^2 e^{-2\alpha\tau} \\ \alpha^2 e^{-\alpha\tau}(e^{-\alpha\tau} a_{23} b + e^{-3\alpha\tau} a_{23} g + 1) \\ \alpha^2 e^{-\alpha\tau}(b + 3e^{-2\alpha\tau} g) \end{pmatrix}$$

$$L_3 = \begin{pmatrix} \alpha^3 e^{-2\alpha\tau}(e^{-\alpha\tau} a_{23} b + e^{-3\alpha\tau} a_{23} g + 3) \\ \alpha^3 e^{-\alpha\tau}(e^{-2\alpha\tau} a_{21} + 3e^{-\alpha\tau} a_{23} b + 7e^{-3\alpha\tau} g + 1) \\ \alpha^3 e^{-\alpha\tau}(b + 9e^{-2\alpha\tau} g) \end{pmatrix}$$

After introducing $\bar{c} = Sy$, $y = (y_1, y_2, y_3)^T$ we get

$$\frac{dy_1}{d\tau} = \varphi_1(\tau)y_3 + d,$$
$$\frac{dy_2}{d\tau} = y_1 + \varphi_2(\tau)y_3, \quad (6.7)$$
$$\frac{dy_3}{d\tau} = y_2 + \varphi_3(\tau)y_3.$$

The change $y_3 = \psi(\tau)$ reduces (6.7) to a linear equation of the third order

$$\psi^{(3)} - \varphi_3(\tau)\psi^{(2)} - \left(2\frac{d\varphi_3}{d\tau} + \varphi_2(\tau)\right)\psi^{(1)}$$
$$- \left(\frac{d^2\varphi_3}{d\tau^2} + \frac{d\varphi_3}{d\tau} + \varphi_1(\tau)\right)\psi = d. \quad (6.8)$$

The variables $y_1$, $y_2$ and $y_3$ are connected to $\psi^{(2)}$, $\psi^{(1)}$ and $\psi$ with

$$y = T\Psi, \Psi = \left(\psi^{(2)}, \psi^{(1)}, \psi\right)^T,$$
$$T = \begin{Vmatrix} 1 & -\varphi_3 & -\left(\frac{d\varphi_3}{d\tau} + \varphi_3\right) \\ 0 & 1 & -\varphi_3 \\ 0 & 0 & 1 \end{Vmatrix}.$$

Let

$$d = \psi^{(3)} - (6 + \varphi_3(\tau))\psi^{(2)}$$
$$- \left(11 + 2\frac{d\varphi_3}{d\tau} + \varphi_2(\tau)\right)\psi^{(1)}$$
$$- \left(6 + \frac{d^2\varphi_3}{d\tau^2} + \frac{d\varphi_3}{d\tau} + \varphi_1(\tau)\right)\psi.$$

Substituting $d(\tau)$ to (6.8) we get the equation with $-1$, $-2$ and $-3$ as roots of the characteristic polynomial. The return to the original variables gives

$$d = \Gamma(\tau)T^{-1}(\tau)S^{-1}(\tau)\bar{c},$$
$$\Gamma = \left(-(6 + \varphi_3), -\left(11 + 2\frac{d\varphi_3}{d\tau} + \varphi_2\right), \right. \quad (6.9)$$
$$\left. -\left(6 + \frac{d^2\varphi_3}{d\tau^2} + \frac{d\varphi_2}{d\tau} + \varphi_1\right)\right).$$

It is obvious that (6.9) provides an exponential decrease of (6.6) solutions. At the final stage we solve the initial value problem for the system obtained from (6.3) after changing the phase coordinates according to (6.5) with the control (6.9). Then we return to the original variables. The initial values for the Cauchy problem are

$$v_1(0) = -\bar{x}_1 - \frac{1}{2}g(\bar{x}_1) - \frac{1}{24}\frac{\partial g}{\partial c_1}(\bar{x}_1)g(\bar{x}_1) - \bar{\bar{g}},$$
$$u_2(0) = g(\bar{x}_1) + \frac{1}{6}\frac{\partial g}{\partial c_1}(\bar{x}_1)g(\bar{x}_1) + \bar{g} + \bar{\bar{g}},$$
$$c_3(0) = 0.$$

During the numerical simulation we have solved the auxiliary system of ODEs constructed from (6.3) and (6.9) after changing the phase coordinates according to (4.5) with the initial values $v_1(0)$, $u_2(0)$, $c_3(0)$ for $\bar{x}_1 = 10$ meters, $r_0 = 7 \cdot 10^6$ meters, $\alpha = 0.1$ and time span $[0; 12.5]$.

Finally we've switched from the auxiliary to the original variables. The plots show the control $u(t)$ and the corresponding transient plots of the phase coordinates $x_1(t)$, $x_2(t)$ and $x(t)$ in respect of the original independent variable $t$.

## 7. Conclusions

The analysis of the Theorem proof shows that the most difficult and time-consuming part of the algorithm implementation can be proceeded with analytical methods of computer algebra packages. The results of the numerical simulation of interorbital flight convince that the method can be used for construction and simulation of various technical objects control systems.

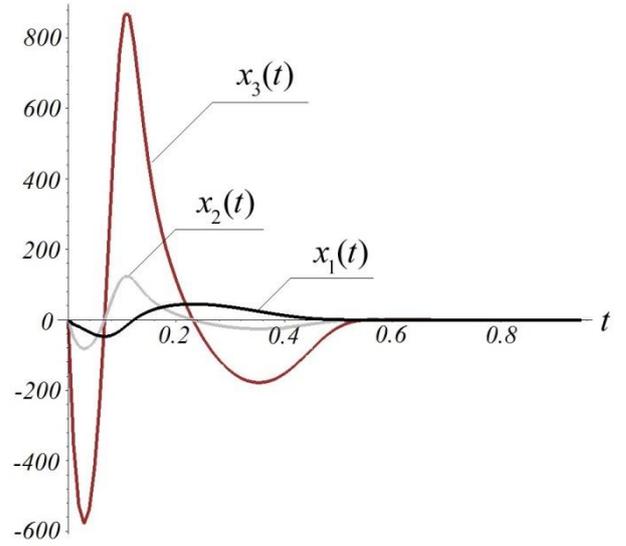

FIGURE 1: Time change of the phase variables

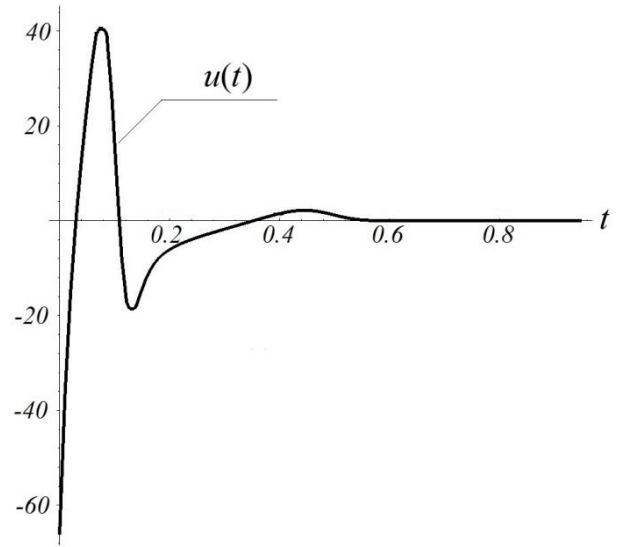

FIGURE 2: Time change of the control $u(t)$